\documentclass[12pt]{article}
\usepackage{natbib,bm}
\usepackage{amsmath}
\usepackage{amssymb}
\usepackage{graphicx}
\usepackage[dvips]{epsfig}
\usepackage{lscape}

\renewcommand{\baselinestretch}{1.1}
\tiny\normalsize
\newcommand{\nmathbf}{\bm}

\def\bfA{\nmathbf A}
\def\bfB{\nmathbf B}
\def\bfC{\nmathbf C}

\def\bfI{\nmathbf I}

\def\bfP{\nmathbf P}

\def\bfR{\nmathbf R}

\def\bfV{\nmathbf V}

\def\bfh{\nmathbf h}

\def\bfs{\nmathbf s}

\def\bfu{\nmathbf u}

\def\bfy{\nmathbf y}
\def\bfz{\nmathbf z}

\def\bftheta  {\nmathbf \theta}

\def\bfmu     {\nmathbf \mu}

\def\bfphi    {\nmathbf \phi}

\def\bfSigma  {\nmathbf \Sigma}

\def\boldfacefake#1{\kern-4pt
   \hbox{ \mathsurround=0pt
   \hbox to 0.4pt{$#1$\hss}\hbox to 0.4pt{$#1$\hss}\hbox {$#1$}}}

\newcommand{\ok}{\hfill\fbox{}}

\newcommand{\btable}{\begin{table}[h]\centering}
\newcommand{\etable}{\end{table}}
\newcommand{\bt}{\begin{parag}\small \let\b=\nsb \let\sb=\nssb \begin{tabular}}
\newcommand{\et}{\end{tabular}\let\b=\nb \let\sb=\nsb\end{parag}}

\newenvironment{parag}{\par}{\par}
\newenvironment{dif}
    {\begin{parag}\small \let\b=\nsb \let\sb=\nssb \begin{parag}}
    {\let\b=\nb \let\sb=\nsb \end{parag}\end{parag}}

\newenvironment{proof}{\begin{dif} \noindent{\em Proof.~}}
            {\ok\vspace*{10pt}\end{dif}}
\newenvironment{exa}{\begin{list}{}
           {\setlength{\leftmargin}{10pt}
            \setlength{\rightmargin}{\leftmargin}}
           \item\begin{ex}\em}{\end{ex}\end{list}}

\newcommand{\be}{\begin{eqnarray}}
\newcommand{\ee}{\end{eqnarray}}
\newcommand{\ba}{\begin{eqnarray*}}
\newcommand{\ea}{\end{eqnarray*}}

\newtheorem{theorem0}{Theorem}
\newtheorem{lemma0}{Lemma}
\newtheorem{remark0}{Remark}
\newtheorem{fact0}{Fact}
\newtheorem{example0}{Example}
\newtheorem{definition0}{Definition}
\newtheorem{corollary0}{Corollary}
\newtheorem{proposition0}{Proposition}
\newtheorem{algorithmY}{Algorithm}

\newcommand{\reals}{\mbox{\rm I\kern-.20em R}}
\newcommand{\sreals}{\mbox{\small \rm I\kern-.20em R}}

\newcommand{\goto}{\rightarrow}

\newcommand{\bdfn}{\begin{dfn}}
\newcommand{\edfn}{\end{dfn}}
\newcommand{\bteo}{\begin{teo}}
\newcommand{\eteo}{\end{teo}}
\newcommand{\bexa}{\begin{exa}}
\newcommand{\eexa}{\end{exa}}
\newcommand{\bdif}{\begin{dif}}
\newcommand{\edif}{\end{dif}}
\newcommand{\bpro}{\begin{proof}}
\newcommand{\epro}{\end{proof}}

\makeatletter
\def\vec{\mathop{\operator@font vec}\nolimits}
\makeatother

\DeclareMathOperator{\diag}{diag}

\begin{document}

\begin{center}
 {\large \bf The Influence of Misspecified Covariance on False Discovery Control when Using Posterior Probabilities}

\bigskip
 Ye Liang$^1$, Joshua D. Habiger$^1$ and Xiaoyi Min$^2$ \\
 {$^1$ \small \it Department of Statistics, Oklahoma State University, Stillwater, OK, USA} \\
 {$^2$ \small \it Department of Mathematics and Statistics, Georgia State University, Atlanta, GA, USA}
\end{center}
\medskip

\begin{abstract}
This paper focuses on the influence of a misspecified covariance structure on false discovery rate for the large scale multiple testing problem. 
Specifically, we evaluate the influence on the marginal distribution of local fdr statistics, which are used in many multiple testing procedures
and related to Bayesian posterior probabilities. Explicit forms of the marginal distributions under both correctly specified and incorrectly specified
models are derived. The Kullback-Leibler divergence is used to quantify the influence caused by a misspecification. Several numerical examples 
are provided to illustrate the influence. A real spatio-temporal data on soil humidity is discussed. 

\end{abstract}

\section{Introduction}
Large scale multiple testing arises from many practical problems, from genetic studies to public health surveillance. \cite{BH1995} introduced
the concept of the false discovery rate (FDR) and proposed a powerful testing procedure, usually referred as the BH procedure. 
The BH procedure relies on a positive dependence assumption 
\citep{BY2001}, while adaptive BH procedures \citep{Storey2004, Liang2012} rely on an independence or weak dependence structure.  
\cite{Efron2007} noted that correlation may result in overly liberal or overly conservative testing procedures. 
Though the BH procedure is valid under different dependence assumptions \citep{Farcomeni2007, Wu2008}, \cite{SunW2009} showed that
failing to model dependence can result in inefficiency. To address that problem, \cite{SunW2009} and \cite{SunW2015}
propose a procedure using local significance index, which is a Bayesian posterior probability. \cite{Efron2001} described the connection between 
the FDR and Bayes procedures, where a posterior probability is referred as the local false discovery rate (Lfdr). \cite{SunW2009}'s local significance index
reduces to the Lfdr under independence. In fact, there is a rich history of using Bayesian approaches for multiplicity adjustment. 
\cite{Scott2006} and \cite{Scott2010} discussed Bayesian multiplicity adjustment in variable selections. \cite{Muller2006} had a comprehensive discussion 
on the connection between the FDR, Bayesian multiple testing, and procedures using posterior probabilities. 

In \cite{SunW2015}, the procedure depending on unknown parameters is called the ``oracle'' procedure, in the sense that we
know all nuisance parameters as an oracle. Although the oracle procedure is proved to control FDR 
at the nominal level and be optimal in terms of false non-discovery rate (FNR), a data driven ``adaptive'' procedure relies on 
correctly specifying the model, including the prior specification, and/or consistent parameter estimation. For large-scale data, dependence 
or covariance, if in Gaussian models, is often estimated based on a structured model. The model choice or the structure choice itself may be 
debatable and parameter estimation remains challenging. For example, in spatial modeling, the estimation of covariance relies on 
structured covariance specifications and in practice, one may have multiple choices of specifications. Intuitively, the choice of specification
will influence the data driven procedure and may eventually lead to different decisions. 

In this paper, we explore the influence of a misspecified covariance structure on the testing procedure.
Specifically, we study the sampling distributions of the local fdr statistics under both correctly and incorrectly specified covariance structures.  
We derive explicit expressions for those distributions under a general model setting. We propose to use the Kullback-Leibler divergence
as a quantitative measure for the influence. We show in both a simulation study and a real application that the influence of 
a misspecification leads to unappealing results.   
The paper is organized as follows. Section 2.1 gives a basic setup of this problem. Section 2.2 gives sampling distributions of the test statistics. 
Section 2.3 provides formulas for computing the Kullback-Leibler divergence. Section 2.4 shows several numerical examples. 
Section 3 provides real data that arose from the Oklahoma soil monitoring network. Section 4 is a discussion. 

\section{Main Results}
\subsection{The Setup}
In this paper, we consider a general model for $m$ observations $\bfy=(y_1,\ldots,y_m)'$:
\be \label{model}
	y_i=\theta_i+\epsilon_i, ~~~ i=1,\ldots,m,
\ee
where $\theta_i$ is the latent state and $\epsilon_i$ is the noise term which independently follows $\mbox{N}(0,\sigma^2)$. The dependence of observations 
is introduced through their latent states $\bftheta=(\theta_1,\ldots,\theta_m)'$. For instance, $\theta_i$ can be a realized spatial process $\theta(s)$, 
for which a spatial dependence structure can be specified. Consider a one-sided hypothesis 
\be \label{hypo}
	H_{0i}: \theta_i \geq \theta_{0i} ~~~\mbox{versus}~~~ H_{1i}:\theta_i < \theta_{0i},
\ee
for every $i$ simultaneously. This type of one-sided hypothesis is often of interest in many practices. 
In spatial epidemiology, one may want to determine which regions or locations have disease rates higher than some given threshold $\theta_0$. 
Here, data will be spatially correlated and each hypothesis will be one-sided. Similarly, in agricultural studies, one may want to determine 
which locations or time periods have soil moisture levels lower than a given threshold $\theta_0$, indicating a risk of drought, 
and the data will be either spatially or temporally correlated and each hypothesis will be one-sided. 
A more general hypothesis would be $H_{0i}: \theta_i \in \Theta_{0i}$ versus $H_{1i}:\theta_i \in \Theta_{1i}$. We do not
consider a precise (or two-sided) hypothesis in this paper, but some comments are given in the discussion section.    

The dependence of latent states $\bftheta$ is usually specified through a prior model. For example, consider a normal-inverse-gamma prior 
on $(\bftheta,\sigma^2)$,
\be \label{prior}
	\bftheta \mid \sigma^2 \sim \mbox{N}_m(\bftheta_0, g\sigma^2\bfSigma) ~~~\mbox{and}~~~
	\sigma^2 \sim \mbox{IG}(\alpha,\beta).
\ee
For simplicity, we assume that $\bfSigma$ is a known covariance structure and $g$ is a known scale parameter.
The use of $g$ here is the same as that in Zellner's $g$-prior for Bayesian variable selection problems. The $g$ value could be fixed, estimated
or have a hyperprior \citep{Liang2008}. The prior specification (\ref{prior}) in fact induces a marginal probability for each hypothesis: 
$P(H_{0i})=P(\theta_{i}\geq \theta_{0i})=0.5$. 

\cite{SunW2009} and \cite{SunW2015} showed that, to control FDR when data are dependent, the posterior probability $h_i=P(H_{0i} \mid \bfy)$
is useful. The posterior probability $h_i$ is viewed as a test statistic, called local index of significance in their work. 
The oracle procedure orders $\bfh^m=(h_{(1)},\ldots,h_{(m)})$ and rejects all $H_{(i)}, i=1,\ldots,k$ such that
 \be \label{control}
 	k=\mbox{max}\left\{ i: \frac{1}{i}\sum_{j=1}^i h_{(j)} \leq \alpha^* \right\},
 \ee 
where $\alpha^*$ is the nominal level. The procedure mimics the Benjamini-Hochberg procedure, in which $p$-value is the test statistic. 
\cite{SunW2009} showed that this oracle procedure controls FDR at level $\alpha^*$ and has the smallest FNR among all
FDR procedures at $\alpha^*$ for a hidden Markov model. \cite{SunW2015} further showed that in a spatial random field model, 
this oracle procedure controls FDR at level $\alpha^*$ and has the smallest missed discovery rate (MDR). 
A data-driven procedure, however, depends on the estimation of other nuisance parameters. The covariance $\bfSigma$ is especially
important in this case as it describes the dependence. Our objective is to determine if the procedure is sensitive when the covariance
is incorrectly specified or estimated, and if so, to quantify the sensitivity.  

\subsection{Sampling Distribution of Test Statistics}
\subsubsection{Known variance of noise}
We now focus on how the distribution of test statistics $(h_1,\ldots,h_m)$ is influenced by a misspecified covariance structure. 
Assume that the data are generated from the true underlying process:
\be \label{true}
	y_i \sim \mbox{N}(\theta_i, \sigma_0^2) ~~~\mbox{and}~~~ \bftheta \sim \mbox{N}(\bftheta_0, \bfSigma_1),
\ee 
for $i=1,\ldots,m$. Assume known $\sigma_0^2$ and consider model (\ref{model}) with priors $\bftheta \sim \mbox{N}(\bftheta_0, g\bfSigma_1)$
and $\bftheta \sim \mbox{N}(\bftheta_0, g\bfSigma_2)$, the latter of which has a misspecified covariance structure. The scale $g$ determines 
the strength of the prior. The intuition is that both $g$ and $\bfSigma_2$ will influence the test statistics $(h_1,\ldots,h_m)$ and FDR control. 

\begin{lemma0} \label{lemma1}
Under the correct covariance $\bfSigma_1$, $h_i$ marginally has the following CDF: 
\be \label{h_dist}
	F(h_i)=\Phi\left[\sqrt{\frac{a_{ii}}{b_{ii}}} \Phi^{-1}(h_i)\right],
\ee
where $a_{ii}$ is the $i$th diagonal element in $\bfA=(1/\sigma_0^2\bfI+1/g\bfSigma_1^{-1})^{-1}$ and 
$b_{ii}$ is the $i$th diagonal element in $\bfB=(\bfI+\sigma_0^2/g\bfSigma_1^{-1})^{-1}(\sigma_0^2\bfI+\bfSigma_1)(\bfI+\sigma_0^2/g\bfSigma_1^{-1})^{-1}$. \\
Under the misspecified covariance $\bfSigma_2$, $\bfA=(1/\sigma_0^2\bfI+1/g\bfSigma_2^{-1})^{-1}$ and
$\bfB=(\bfI+\sigma_0^2/g\bfSigma_2^{-1})^{-1}(\sigma_0^2\bfI+\bfSigma_1)(\bfI+\sigma_0^2/g\bfSigma_2^{-1})^{-1}$.
\end{lemma0}  

\begin{figure}[ht]
\begin{center}
\includegraphics[scale=0.65,angle=270]{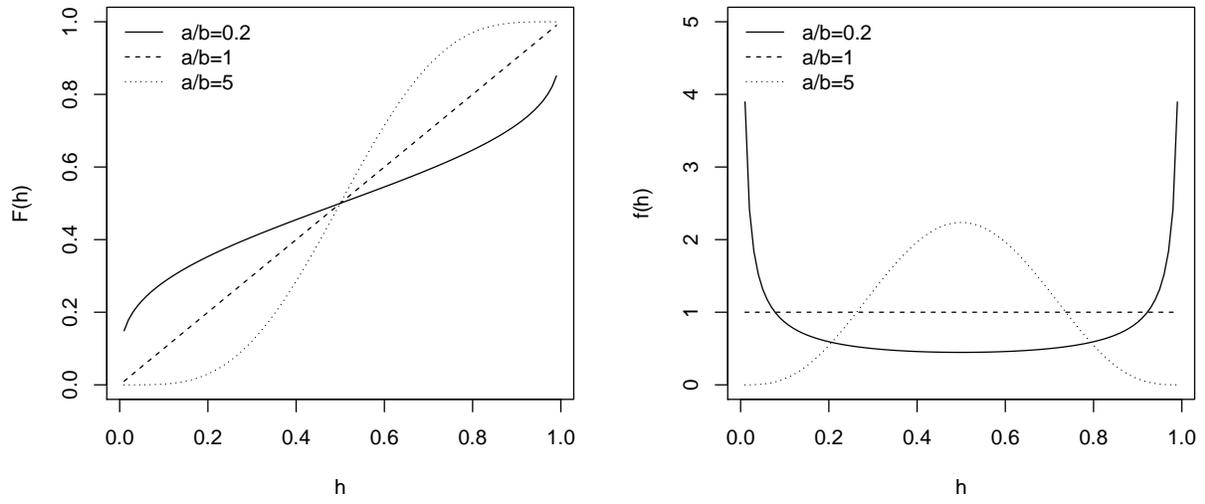}
\caption{Marginal CDF and pdf for $h_i$.}
\label{fig:marginal}
\end{center}
\end{figure}

Lemma \ref{lemma1} shows explicitly how the sampling distribution is altered by a misspecification. 
Observe that $F(h_i)$ is completely determined by the ratio $a_{ii}/b_{ii}$. Figure \ref{fig:marginal} shows different shapes of both CDF
and pdf under different ratio values. Note that when $a_{ii}=b_{ii}$, the sampling distribution is Uniform$(0,1)$. Consider $g\goto \infty$, for which
the prior becomes noninformative, then $a_{ii}/b_{ii} \goto \sigma_0^2/(\sigma_0^2+\sigma^2_{1,ii})$, where $\sigma^2_{1,ii}$ is the $i$th
diagonal element in $\bfSigma_1$. The ratio becomes irrelevant to the misspecified covariance $\bfSigma_2$. 
In other words, the behaviors of the correct specification and the incorrect specification will be similar when $g$ is large, which seems intuitive. 
Because FDR procedures based on $h_i$ reject $H_{0i}$ if $h_i\leq C$ for some $C$, it is important to note that $P(H_i \leq h_i)=F(h_i)$
can be substantially influenced by covariance misspecification. For example, in Figure \ref{fig:marginal}, we see $F(0.2)$ ranges from $0.03$
to $0.35$. It should be noted that some recommend rejecting $H_{0i}$ if $h_i \leq 0.2$ \citep{Efron2012}. 
In this setting, the rejection probability ranges from $0.03$ to $0.35$.   

Unlike the independent case, $h_1,\ldots,h_m$ are now dependent when the data are dependent. The change under a misspecified structure 
is revealed in their joint distribution.

\begin{theorem0} \label{thm1}
Using the same definition for $\bfA$ and $\bfB$ in Lemma \ref{lemma1}, $h_1,\ldots,h_m$ have a joint CDF 
\ba \label{h_joint1}
	F(h_1,\ldots,h_m)=\Phi_b^m\left[\sqrt{\frac{a_{11}}{b_{11}}} \Phi^{-1}(h_1),\cdots,\sqrt{\frac{a_{mm}}{b_{mm}}} \Phi^{-1}(h_m)\right], 
\ea
where $\Phi_b^m$ is the CDF for a multivariate normal $\mbox{N}_m({\bf0}, \bfP_b)$, and $\bfP_b$ is the correlation matrix of $\bfB$. 
\end{theorem0}

The joint distribution of $h_1,\ldots,h_m$ represents a multivariate surface in the space $[0,1]^m$. Notice that jointly not only the ratio $a_{ii}/b_{ii}$
plays a role but also the correlation structure of $\bfB$. Under a misspecified covariance structure, $\bfP_b$ will be altered as well. 
However, still, as $g\goto \infty$, $\bfB \goto \sigma_0^2\bfI+\bfSigma_1$, and there is no misspecification effect.

\subsubsection{Unknown variance of noise}  
Assume the underlying data generating process (\ref{true}). Suppose $\sigma_0^2$ is unknown and consider specifying model (\ref{model})
with prior (\ref{prior}). As before, a correct covariance structure is $\bfSigma_1$ and a misspecified structure is $\bfSigma_2$. 

 \begin{theorem0} \label{thm2}
Under the correct covariance $\bfSigma_1$, the test statistics $h_1, \ldots, h_m$ jointly have 
the following CDF: 
 \be \label{h_joint2}
 	F(h_1, \ldots, h_m)=\Xi^m_{a, b}\left[ \sqrt{\frac{a_{11}}{b_{11}}} \Psi^{-1}_{m+2\alpha}(h_1), \cdots, 
	\sqrt{\frac{a_{mm}}{b_{mm}}} \Psi^{-1}_{m+2\alpha}(h_m) \right],
 \ee
 where $\bfA=(\bfI+1/g\bfSigma_1^{-1})^{-1}$, $\bfB=(\bfI+1/g\bfSigma_1^{-1})^{-1}(\sigma_0^2\bfI+\bfSigma_1)(\bfI+1/g\bfSigma_1^{-1})^{-1}$, 
 $\Psi_{m+2\alpha}$ is the CDF for a univariate $t$-distribution with degrees of freedom $m+2\alpha$,  and $\Xi^m_{a, b}$ 
 is the CDF for the following random vector
 \ba
 	\sqrt{\frac{m+2\alpha}{\bfz_{b}' \bfC \bfz_{b} + 2\beta}} ~\bfz_{b},
 \ea
 where $\bfz_{b} \sim \mbox{N}({\bf0}, \bfP_b)$ and $\bfC=(\diag\bfB)^{-1/2}(\bfA^{-2}-\bfA^{-1})(\diag\bfB)^{-1/2}$, 
 where $\diag$ here denotes a diagonal matrix. \\
 Under the misspecified covariance $\bfSigma_2$, $\bfA=(\bfI+1/g\bfSigma_2^{-1})^{-1}$
 and $\bfB=(\bfI+1/g\bfSigma_2^{-1})^{-1}(\sigma_0^2\bfI+\bfSigma_1)(\bfI+1/g\bfSigma_2^{-1})^{-1}$. 
\end{theorem0}

The main difference between Theorem \ref{thm2} and Theorem \ref{thm1} is that $\Xi^m_{a,b}$ has a more complicated form than $\Phi^m_b$.
The ratio $a_{ii}/b_{ii}$ still plays a role in the joint distribution and $\Xi^m_{a,b}$ will be affected by a misspecified structure. 

\subsection{Kullback-Leibler Divergence}
Denote the sampling distribution of $\bfh$ under the correct covariance specification as $f_{\mbox{\scriptsize cor}}(\bfh)$ and 
under the misspecified covariance as $f_{\mbox{\scriptsize mis}}(\bfh)$. We evaluate the influence of the misspecification by 
the Kullback-Leibler (KL) divergence
\be \label{KLDiv}
	D_{\mbox{\scriptsize KL}}(f_{\mbox{\scriptsize cor}} ~\|~  f_{\mbox{\scriptsize mis}})
	= \int f_{\mbox{\scriptsize cor}}(\bfh) \log \frac{f_{\mbox{\scriptsize cor}}(\bfh)}{f_{\mbox{\scriptsize mis}}(\bfh)} d\bfh.
\ee
The KL divergence here can be interpreted as the information loss when using $f_{\mbox{\scriptsize mis}}$ to approximate 
$f_{\mbox{\scriptsize cor}}$. The following two corollaries are useful to approximate the KL divergence in our cases. 

\begin{corollary0} \label{coro1}
As a consequence of Lemma \ref{lemma1}, the marginal density function of $h_i$ is 
\ba
	f(h_i)=\sqrt{r_i}\exp \left\{\frac{1}{2}(1-r_i)\phi_i^2 \right\},
\ea
where $r_i=a_{ii}/b_{ii}$ and $\phi_i=\Phi^{-1}(h_i)$.
\end{corollary0}

\begin{corollary0} \label{coro2}
As a consequence of Theorem \ref{thm1}, the joint density function of $(h_1,\ldots,h_m)'$ is
\ba
	f(h_1,\ldots,h_m)=\left| \bfR^{\frac{1}{2}}\bfP_b^{-1} \bfR^{\frac{1}{2}} \right|^{\frac{1}{2}} 
	\exp \left\{ \frac{1}{2}\bfphi'(\bfI-\bfR^{\frac{1}{2}}\bfP_b^{-1} \bfR^{\frac{1}{2}})\bfphi \right\},
\ea
where $\bfphi=(\phi_1,\ldots,\phi_m)'=(\Phi^{-1}(h_1),\ldots,\Phi^{-1}(h_m))'$ and $\bfR=\diag(r_1,\ldots,r_m)$. 
\end{corollary0}

With Corollary \ref{coro2}, the term $\log \{f_{\mbox{\scriptsize cor}}(\bfh) / f_{\mbox{\scriptsize mis}}(\bfh)\}$ in expression (\ref{KLDiv}) is analytically available. 
Thus, the KL divergence $D_{\mbox{\scriptsize KL}}$ can be easily evaluated using Monte Carlo approximation. Notice that we can draw from 
$f_{\mbox{\scriptsize cor}}$ exactly given the underlying model because $h_i=P(H_{0i} \mid \bfy)$ is analytically available.
To draw a sample of $\bfh$ from $f_{\mbox{\scriptsize cor}}$, draw $\bfy$ from the underlying model and then calculate $\bfh$. Suppose we obtain
a sample $\bfh_1,\ldots,\bfh_L$, the approximation is
\ba
	D_{\mbox{\scriptsize KL}}(f_{\mbox{\scriptsize cor}} ~\|~  f_{\mbox{\scriptsize mis}}) \approx
	\frac{1}{L}\sum_{l=1}^L \log \frac{f_{\mbox{\scriptsize cor}}(\bfh_l)}{f_{\mbox{\scriptsize mis}}(\bfh_l)}. 
\ea
A misspecified covariance will change the defined matrices $\bfA$ and $\bfB$, and consequently change matrices $\bfR$ and $\bfP_b$. 
Notice the relationship $\bfR^{1/2}\bfP_b^{-1} \bfR^{1/2}=(\diag^{1/2}\bfA)\bfB^{-1} (\diag^{1/2}\bfA)$. Since the KL divergence will generally
increase as the dimension $m$ increases (in the independent case it is simply a sum of individual dimensions), we may also consider
a relative measure of influence $D_{\mbox{\scriptsize KL}}/m$. 

The KL divergence is computable under the general model (\ref{model}) with $\bfSigma_1$ and $\bfSigma_2$ provided.
This easy-to-compute measure can be used to quantify the influence of a misspecification. In practice, when there are multiple
candidate covariances, we may assess the KL divergence between those candidates.  

\subsection{Numerical Examples}
We now consider two examples of misspecified covariance. In each example, without loss of generality, we set $m=900$ and $\sigma_0^2=0.25$. 
We will numerically evaluate the KL divergence and perform a simulation study, in which we estimate FDR and FNR with Monte Carlo replications of $1000$. 

\begin{example0}
Positive spatial covariance vs. Independence. \\
Consider a regular spatial grid with unit distance one for generating the latent states $\bftheta$. The true covariance $\bfSigma_1=\{\sigma^2_{1,ij}\}$ has a positive
decaying structure determined by an exponential covariance function $\sigma^2_{1,ij}=\exp\{-\|s_i-s_j\|/\rho\}$ with $\rho=5$, where $s$ represents a location. 
A misspecified covariance is $\bfSigma_2=\bfI$. To get a rough idea of $r_i=a_{ii}/b_{ii}$, let g=1 and compute $\bfA$ and $\bfB$. Under the correct specification, 
$r_i$ ranges from $0.12$ to $0.16$, and under the misspecification, $r_i=0.25$. 
\end{example0}

Note that this misspecification is essentially to ignore the dependence and treat data as independent observations. 
This is quite common in practice, where domain scientists are hesitant to model complex covariance structures, though evidences
suggest that data may be correlated. 

\begin{example0}
Negative AR(2) covariance vs. Positive AR(2) covariance. \\
Consider a time series for generating the latent states $\bftheta$. The true covariance $\bfSigma_1$ is determined by an AR(2) process: 
$\theta_i=\rho_1 \theta_{i-1}+\rho_2 \theta_{i-2}+\varepsilon_i, ~\varepsilon_i \sim \mbox{N}(0,1)$, and $\rho_1=1.5$ and $\rho_2=-0.9$. The autocorrelation
function of this specification has an oscillating pattern (mixed positive and negative values in $\bfSigma_1$). A misspecified covariance $\bfSigma_2$
is chosen to be the covariance for an AR(2) process with $\rho_1=0.6$ and $\rho_2=0.3$, whose autocorrelation is always positive. 
To get a rough idea of $r_i=a_{ii}/b_{ii}$, let g=1 and compute $\bfA$ and $\bfB$. Under the correct specification, $r_i$
ranges from $0.088$ to $0.14$, and under the misspecification, $r_i$ ranges from $0.20$ to $0.25$. 
\end{example0}

Note that this example shows a scenario where model is correctly specified but parameter estimates are wrong. 
This example also compares a covariance matrix containing negative values with a covariance matrix containing all positive values. 

\begin{figure}[t]
\begin{center}
\includegraphics[scale=0.7,angle=270]{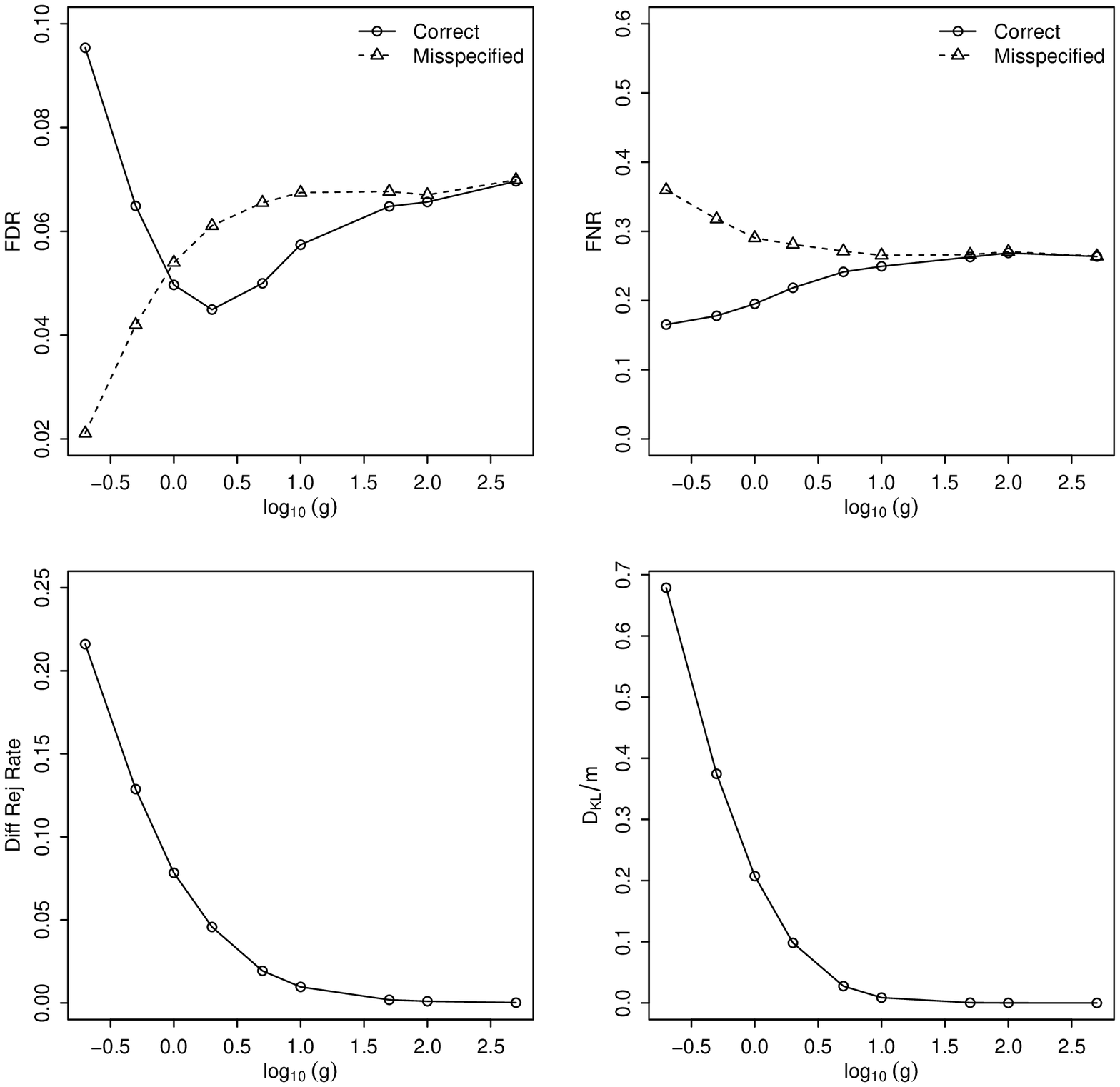}
\caption{Results for Example 1. FDR, FNR, Difference of rejection rate ($\mbox{rate}_{\mbox{\scriptsize cor}}-\mbox{rate}_{\mbox{\scriptsize mis}}$) 
and $D_{\mbox{\scriptsize KL}}/m$. The sequence of $g$ is from $0.2$ to $500$. The nominal level is $0.05$.}
\label{fig:ex1}
\end{center}
\end{figure}

\begin{figure}[t]
\begin{center}
\includegraphics[scale=0.7,angle=270]{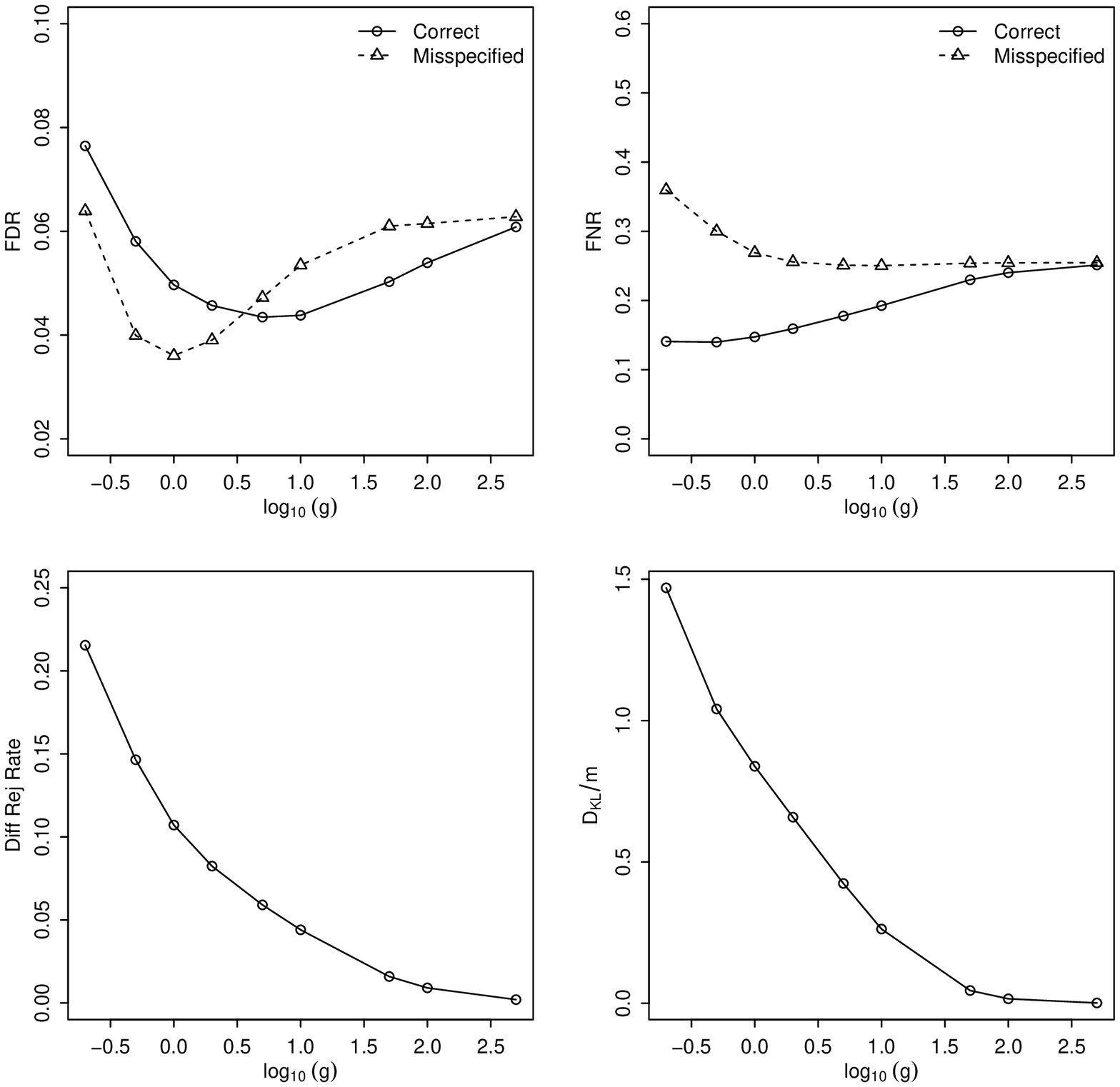}
\caption{Results for Example 2. FDR, FNR, Difference of rejection rate ($\mbox{rate}_{\mbox{\scriptsize cor}}-\mbox{rate}_{\mbox{\scriptsize mis}}$) 
and $D_{\mbox{\scriptsize KL}}/m$. The sequence of $g$ is from $0.2$ to $500$. The nominal level is $0.05$.}
\label{fig:ex2}
\end{center}
\end{figure}

Results of Example 1 and 2 are shown in Figure \ref{fig:ex1} and Figure \ref{fig:ex2}. We choose different $g$ values representing different strengths of
information brought in by the prior dependence. Four plots are shown in each result: the estimated FDR, the estimated FNR, the difference between the 
rejection rate of the correct specification and that of the misspecification, i.e. $(\# discovery_{\mbox{\scriptsize cor}} - \# discovery_{\mbox{\scriptsize mis}})/m$, 
and $D_{\mbox{\scriptsize KL}}/m$. We can reach the following conclusions 
from these plots. First, as $g$ increases, both the KL divergence and the difference between rejection rates decrease, and also both the FDR 
and the FNR become closer, all of which are as expected, indicating a decreasing misspecification influence. 
Second, the FNR under the misspecification is universally higher than under the correct specification. 
Hence, misspecification results in an inefficient procedure. Also notice that the procedure tends to give less discoveries under the misspecification
than under the correct specification. 
Last but not least, the FDR change is not monotonic and the comparison between the two specifications is profound. Notice that the nominal level 
is $0.05$ and $g=1$ (or $\log_{10}g=0$) represents a ``true'' scale. When both the structure $\bfSigma_1$ and the scale $g$ are correct, 
from the top left plot in both results, we can see that the FDR is controlled at the nominal level. This seems to suggest that
a correctly estimated scale of prior dependence is desired, which should neither be too strong nor too weak. And this correct scale will only
work as expected if the covariance structure is correct as well. 

\begin{example0}
Positive spatial covariance revisit. \\
Revisit Example 1. Consider to fix $g=1$. The correct specification is exactly the same as the underlying model with $\rho=5$. 
For the misspecification, let the spatial range parameter $\rho$ change from $0.1$ to $20$, representing the strength of dependence. 
\end{example0}

\begin{figure}[t]
\begin{center}
\includegraphics[scale=0.8,angle=270]{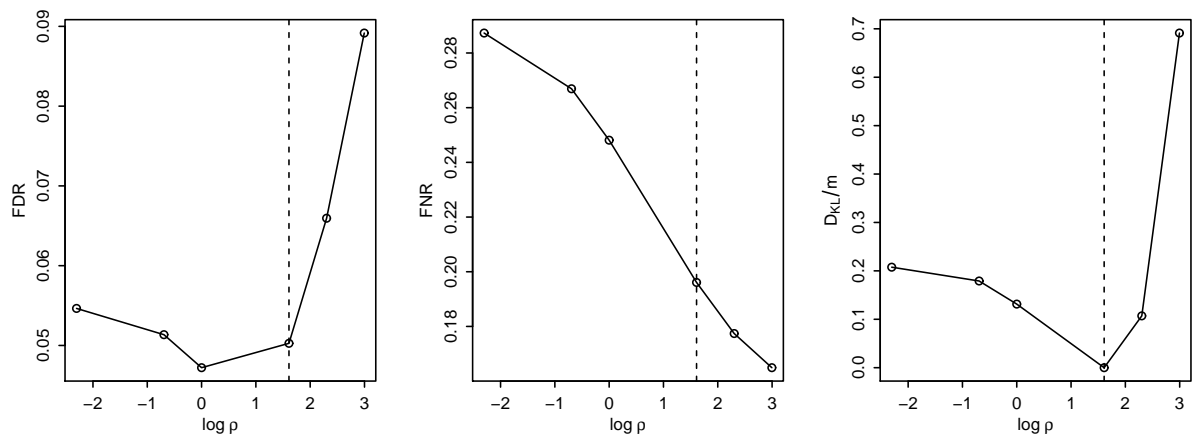}
\caption{Results for Example 3. FDR, FNR and $D_{\mbox{\scriptsize KL}}/m$. The dashed vertical line is the correct value $\rho=5$. 
The sequence of $\rho$ is $(0.1, 0.5, 1, 5, 10, 20)$. The nominal level is $0.05$.}
\label{fig:ex3}
\end{center}
\end{figure}

Results of Example 3 are shown in Figure \ref{fig:ex3}. The FDR is maintained at $0.05$ only when $\rho$ has the correct value.  
The pattern of change in FDR is also reflected in the KL divergence plot. In this example, using $D_{\mbox{\scriptsize KL}}/m$ as a measure of influence 
seems to be reasonable. The FNR, on the other hand, monotonically decreases as the dependence goes stronger. We shall note here that,
when $g$ is large, the KL divergence becomes sensitive to detect a misspecification as the measure approaches zero.
However, in that case, the prior is vague and the influence on FDR is negligible. 

\section{Real Data: Soil Relative Humidity}
Oklahoma Mesonet \citep{Illston2008} is a comprehensive observatory network monitoring environmental variables across the state. One of the focuses
of the network is soil moisture. Extreme weather conditions, especially drought, severely impact Oklahoma's agriculture, which is a leading economy
of the state. Soil moisture is fundamentally important to many hydrological, biological and biogeochemical processes. The information is valuable to a wide 
range of government agencies and private companies. We take a small dataset from their data warehouse as an example of real application. 
Consider only one variable here: the {\it relative humidity}, ranging from $0$ to $100$ percent. The dataset consists of monthly averages in 2014
for 108 monitoring stations, which is in total 1,296 measurements. Consider each hypothesis being $H_{0i}: \theta_i \geq 50$ versus $H_{0i}: \theta_i < 50$ 
for detecting low humidity times and locations. 

We consider two different specifications for the dependence structure. Consider a spatio-temporal process, for a spatial location $\bfs$
and a time point $t$, $y(\bfs;t)=\theta(\bfs;t)+\epsilon(\bfs;t)$,
where $\epsilon(\bfs;t)$ is pure error process with $\mbox{N}(0,\sigma^2)$, and $\theta(\bfs;t)$ is a stationary Gaussian process with a constant mean $\mu$ 
and a separable covariance function:
\ba
	C(h; \tau) = \delta C^{(s)}(h) C^{(t)}(\tau),
\ea
where $h=\|\bfs-\bfs'\|$ and $\tau=|t-t'|$ are both Euclidean distances. Specify $C^{(s)}(h)=\exp\{-h/\rho\}$ and $ C^{(t)}(\tau)=\alpha^{\tau}$.
Specify priors for parameters: $\sigma^2 \sim \mbox{IG}(1,1)$, $\delta \sim \mbox{IG}(1,1)$, $\rho \sim \mbox{Uni}(0,+\infty)$ and 
$\alpha \sim \mbox{Uni}(0,1)$. Posterior distributions are obtained through standard Markov chain Monte Carlo (MCMC). 
We ensure that the chain is long enough to converge and take $10,000$ MCMC samples.  
Posterior probabilities $P(H_{0i}\mid \bfy)$ are approximated with posterior samples.  

The second model is specified that the process is independent over time but the variance is not stationary over time: 
\ba
	C(h; t, t') = 
	\begin{cases}
		\delta_t C^{(s)}(h), & \text{if } t=t'; \\
		0, & \text{otherwise.}
	\end{cases}
\ea
This is a different specification from Model 1 and neither model is a reduced case of the other. 
We use the same specification for $C^{(s)}(h)$ as in Model 1 and also the same prior distributions for $\sigma^2$, $\delta_t$ and $\rho$. 
We take $10,000$ MCMC samples and posterior probabilities are approximated with MCMC samples. 

\begin{figure}[t]
\begin{center}
\includegraphics[scale=0.7,angle=270]{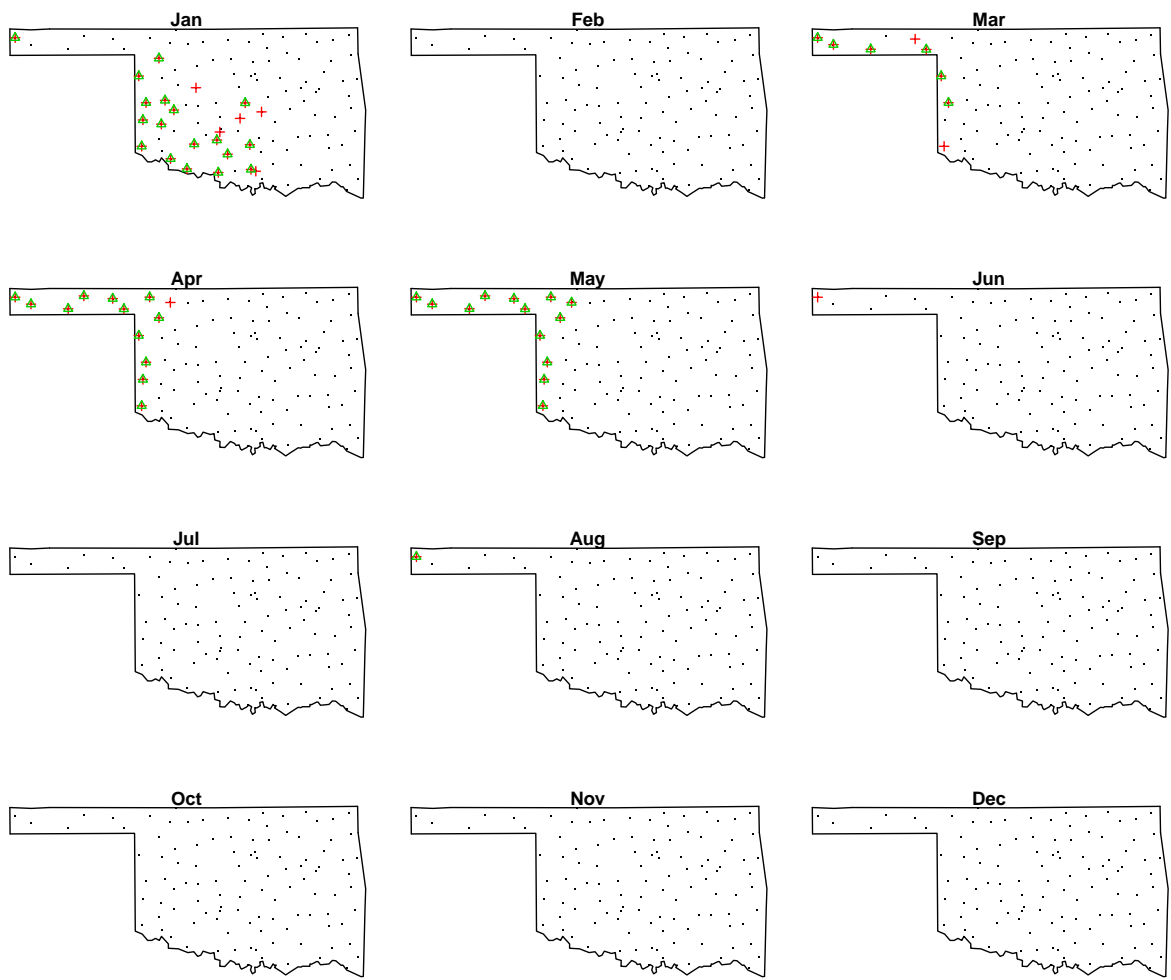}
\caption{Rejected sites for the year of 2014 in Oklahoma, from January to December. Black $(\cdot)$ are the observations,  
	red $(+)$ are rejections using Model 1, and green $(\triangle)$ are rejections using Model 2. }
\label{fig:data}
\end{center}
\end{figure}

For both models, we follow the FDR control procedure (\ref{control}) at a nominal level $0.05$. Results of both are shown in Figure \ref{fig:data}. 
Note that drought impacts Oklahoma mostly in the western areas. Both results seem reasonable and meaningful for practitioners
and they overlap on most decisions. However, we do observe that, at nine time/location points, they do not agree with each other.
Those nine points are all rejected in Model 1 but neither in Model 2. Table \ref{tab:data} shows the observed values and model inferred 
upper credible bounds at the given nominal level for the nine disagreed points. We can see that all nine points are boundary cases
and Model 2 results in higher upper bounds than Model 1, causing the disagreed decisions. 
Such disagreed decisions will likely cause confusions in practice. As one must assume normality in the first place before performing 
a small sample $t$ test, we believe that, in a good practice, it is necessary to clearly assume and carefully check the model specification 
before using posterior probabilities from the model for testing. 

\section{Discussion}
In this paper, we explore the influence of a misspecified covariance structure on the multiple testing procedure using Bayesian posterior probabilities. 
We explicitly show the influence on the test statistics and discuss the KL divergence as a measure of that influence. 
We see from a simulation study that both the correct strength of dependence and structure of dependence are necessary to ensure control of the FDR
 at the nominal level. We also see that misspecified covariance can significantly impact efficiency, in terms of FNR. 
From a real application, we see that different covariance specifications can result in different decisions. 

This paper does not cover any discussion on a precise (or two-sided) hypothesis:
$H_{0i}: \theta_{i}=0$ versus $H_{1i}: \theta_{i}\neq 0$. In that scenario, a mixture model is often assumed: $f(y)=p_0f_0(y)+p_1f_1(y)$. 
The local fdr by \cite{Efron2001} is $p_0f_0(y)/f(y)$ under independence. When data are dependent, it is unclear how to properly incorporate 
the dependence into the mixture model. One practical example given by \cite{Brown2014} specifies: $y_i \sim \mbox{N}(\gamma_i \mu_i, \sigma^2)$, 
$\gamma_i \sim \mbox{Bern}(1-p)$ and $\bfmu \sim \mbox{CAR}(\rho, \tau^2)$, which mimics the independent model in \cite{Scott2006}. 
In the Bayesian framework, to compute $P(H_{0i}\mid \bfy)$, we would need $P(H_{0i})$, $P(H_{1i})$ and priors $\pi_{i0}(\theta_i)$ under $H_{0i}$
and $\pi_{i1}(\theta_i)$ under $H_{1i}$. Moreover, $\pi_{i0}(\theta_i)$ and $\pi_{i1}(\theta_i)$ should have a dependence structure for $i=1,\ldots,m$, 
in some way. A misspecified covariance (or model) would be worth further investigation in this setting. 

\appendix
\section{Appendix}
\subsection{Proof of Lemma \ref{lemma1}}
\begin{proof}
First, given the underlying true model (\ref{true}), it is straightforward to derive the true marginal distribution for $\bfy$:
\be \label{ymar1}
	p(\bfy)=\int_{\bftheta}p(\bfy \mid \bftheta)p(\bftheta) d\bftheta = \mbox{N}(\bftheta_0, \sigma_0^2\bfI+\bfSigma_1).
\ee
If we estimate the posterior using the correct covariance structure, we will have the following posterior distribution: 
$\bftheta \mid \bfy \sim \mbox{N}(\bftheta^{(y)}, \bfSigma^{(y)})$, where
$\bftheta^{(y)}=(1/\sigma_0^2\bfI+1/g\bfSigma_1^{-1})^{-1}(1/\sigma_0^2 \bfy + 1/g\bfSigma_1^{-1}\bftheta_0)$ and 
$\bfSigma^{(y)}=(1/\sigma_0^2\bfI+1/g\bfSigma_1^{-1})^{-1}$. Marginally, $\theta_i \mid \bfy \sim \mbox{N}(\theta_i^{(y)}, \bfSigma_{ii}^{(y)})$.
Then,
\ba
	H_i=P(H_{0i} \mid \bfy)=P(\theta_i \geq \theta_{0i} \mid \bfy)=\Phi\left( \frac{\theta_i^{(y)}-\theta_{0i}}{\sqrt{\Sigma_{ii}^{(y)}}} \right).
\ea
Using (\ref{ymar1}), we have the marginal distribution: $\bftheta^{(y)} \sim \mbox{N}(\bftheta_0, \bfB)$, where 
$\bfB=(\bfI+\sigma_0^2/g\bfSigma_1^{-1})^{-1}(\sigma_0^2\bfI+\bfSigma_1)(\bfI+\sigma_0^2/g\bfSigma_1^{-1})^{-1}$. 
Note that $\bfSigma^{(y)}$ is free of $\bfy$, so let $\bfA=\bfSigma^{(y)}$. 
Marginally, $\theta_i^{(y)} \sim \mbox{N}(\theta_{0i}, b_{ii})$ which leads to
\ba
	 \frac{\theta_i^{(y)}-\theta_{0i}}{\sqrt{\Sigma_{ii}^{(y)}}} \sim \mbox{N}(0, b_{ii}/a_{ii}). 
\ea
Now, under the true covariance, the CDF for $H_i$ is given by
\ba
	F(h_i) &=& P\left[ \Phi\left( \frac{\theta_i^{(y)}-\theta_{0i}}{\sqrt{\Sigma_{ii}^{(y)}}} \right) \leq h_i \right] \\
	&=& P\left[  \frac{\theta_i^{(y)}-\theta_{0i}}{\sqrt{\Sigma_{ii}^{(y)}}} \leq \Phi^{-1}(h_i) \right] \\
	&=& \Phi\left[\sqrt{\frac{a_{ii}}{b_{ii}}} \Phi^{-1}(h_i)\right].
\ea
If a misspecified covariance $\bfSigma_2$ is used to estimate the posterior, then in the posterior distribution 
$\bftheta \mid \bfy \sim \mbox{N}(\bftheta^{(y)}, \bfSigma^{(y)})$, we will have $\bfSigma_2$, instead of $\bfSigma_1$ 
in both $\bftheta^{(y)}$ and $\bfSigma^{(y)}$. As a consequence, $\bfA=(1/\sigma_0^2\bfI+1/g\bfSigma_2^{-1})^{-1}$
and $\bfB=(\bfI+\sigma_0^2/g\bfSigma_2^{-1})^{-1}(\sigma_0^2\bfI+\bfSigma_1)(\bfI+\sigma_0^2/g\bfSigma_2^{-1})^{-1}$. 
The rest remains the same.
\end{proof}

\subsection{Proof of Theorem \ref{thm1}}
\begin{proof}
According to the proof in Lemma \ref{lemma1}, $\bftheta^{(y)} \sim \mbox{N}(\bftheta_0, \bfB)$. Hence, 
$\diag^{-1/2}(\bfB)(\bftheta^{(y)}-\bftheta_0) \sim \mbox{N}({\bf0}, \bfP_b)$, where $\bfP_b$ is the correlation matrix of $\bfB$.
Equivalently, 
\ba
	\left( \frac{\theta_1^{(y)}-\theta_{01}}{\sqrt{b_{11}}}, \cdots, \frac{\theta_m^{(y)}-\theta_{0m}}{\sqrt{b_{mm}}} \right)' \sim \mbox{N}({\bf0}, \bfP_b).
\ea
Then, the joint CDF of $(H_1, \ldots, H_m)'$ is
\ba
	F(h_1,\ldots,h_m) &=& P(H_1 \leq h_1, \ldots, H_m \leq h_m) \\
		&=& P\left[ \Phi\left( \frac{\theta_1^{(y)}-\theta_{01}}{\sqrt{a_{11}}}\right)\leq h_1, \cdots, 
			\Phi\left( \frac{\theta_m^{(y)}-\theta_{0m}}{\sqrt{a_{mm}}}\right)\leq h_m \right] \\
		&=& P\left[ \frac{\theta_1^{(y)}-\theta_{01}}{\sqrt{b_{11}}} \leq \sqrt{\frac{a_{11}}{b_{11}}}\Phi^{-1}(h_1), \cdots, 
			\frac{\theta_m^{(y)}-\theta_{0m}}{\sqrt{b_{mm}}} \leq \sqrt{\frac{a_{mm}}{b_{mm}}}\Phi^{-1}(h_m) \right] \\
		&=& \Phi_b^m\left[\sqrt{\frac{a_{11}}{b_{11}}} \Phi^{-1}(h_1),\cdots,\sqrt{\frac{a_{mm}}{b_{mm}}} \Phi^{-1}(h_m)\right].
\ea
\end{proof}

\subsection{Proof of Theorem \ref{thm2}}
\begin{proof}
 First, given the underlying true process (\ref{true}), marginally $\bfy \sim \mbox{N}(\bftheta_0, \sigma_0^2\bfI+\bfSigma_1)$.
 
 If the correct covariance is used, we have the posterior distribution
 \ba
	p(\bftheta \mid \bfy) &\propto& \int_0^{\infty}p(\bfy \mid \bftheta)p(\bftheta \mid \sigma^2)p(\sigma^2) d\sigma^2 \\
	&\propto& \left[ (\bfy-\bftheta)'(\bfy-\bftheta)+1/g(\bftheta-\bftheta_0)'\bfSigma^{-1}(\bftheta-\bftheta_0) + 2\beta \right]^{-m-\alpha} \\
	&\propto& \left[ (\bftheta-\bftheta^{(y)})'(\bfI+1/g\bfSigma^{-1})(\bftheta-\bftheta^{(y)}) + 
		(\bfy-\bftheta_0)'(\bfI+g\bfSigma)^{-1}(\bfy-\bftheta_0) +2\beta \right]^{-m-\alpha} \\
	&\propto& \left[ 1+\frac{1}{m+2\alpha} (\bftheta-\bftheta^{(y)})'(\bfV^{(y)})^{-1}(\bftheta-\bftheta^{(y)}) \right]^{-\frac{m+(m+2\alpha)}{2}},
\ea
which is $t_m(m+2\alpha, \bftheta^{(y)}, \bfV^{(y)})$, with location 
$\bftheta^{(y)}=\left( \bfI+1/g\bfSigma^{-1} \right)^{-1} \left( \bfy+1/g\bfSigma^{-1}\bftheta_0 \right)$ and scale 
$\bfV^{(y)}=(m+2\alpha)^{-1} \left[2\beta + (\bfy-\bftheta_0)' (\bfI+g\bfSigma)^{-1} (\bfy-\bftheta_0) \right]\left( \bfI+1/g\bfSigma^{-1} \right)^{-1}$.
Similarly as Lemma \ref{lemma1}, define $\bfA=(\bfI+1/g\bfSigma_1^{-1})^{-1}$ and 
$\bfB=(\bfI+1/g\bfSigma_1^{-1})^{-1}(\sigma_0^2\bfI+\bfSigma_1)(\bfI+1/g\bfSigma_1^{-1})^{-1}$. 

Marginally $\theta_i \mid \bfy \sim t_1(m+2\alpha, \theta_i^{(y)}, V_{ii}^{(y)})$, then each test statistic is 
 \ba
 	H_i = P(H_{0i} \mid \bfy) = P(\theta_i \geq \theta_{i0} \mid \bfy) = \Psi_{m+2\alpha}\left( \frac{\theta_i^{(y)}-\theta_{i0}}{\sqrt{V_{ii}^{(y)}}} \right). 
 \ea
 
 In order to find the joint CDF, we need the joint distribution of 
 $((\theta_1^{(y)}-\theta_{10})/\sqrt{V_{11}^{(y)}}, \ldots, (\theta_m^{(y)}-\theta_{m0})/\sqrt{V_{mm}^{(y)}})'$,
 which can be re-written as 
 \be \label{target}
 	\left( \diag \bfV^{(y)} \right)^{-\frac{1}{2}}(\bftheta^{(y)}-\bftheta_0) =  \sqrt{\frac{m+2\alpha}{(\bfy-\bftheta_0)' (\bfI+g\bfSigma_1)^{-1} (\bfy-\bftheta_0)+2\beta}}
	(\diag \bfA)^{-\frac{1}{2}} (\bftheta^{(y)}-\bftheta_0),
 \ee

Given the marginal distribution of $\bfy$, we have $\bftheta^{(y)}-\bftheta_0=\bfA(\bfy-\bftheta_0) \sim \mbox{N}(\bftheta_0, \bfB)$,
or, $\bftheta^{(y)}-\bftheta_0=(\diag\bfB)^{1/2}\bfz_b$. The quadratic term in (\ref{target}) is 
\ba
 	(\bfy-\bftheta_0)' (\bfI+g\bfSigma_1)^{-1} (\bfy-\bftheta_0) &=& (\bfy-\bftheta_0)' (\bfI-\bfA) (\bfy-\bftheta_0) \\
	&=& \bfz_b' (\diag\bfB)^{-1/2}(\bfA^{-2}-\bfA^{-1})(\diag\bfB)^{-1/2} \bfz_b. 
 \ea
The equation (\ref{target}) is then
\ba
	\left( \diag \bfV^{(y)} \right)^{-\frac{1}{2}}(\bftheta^{(y)}-\bftheta_0) &=& \sqrt{\frac{m+2\alpha}{\bfz_{b}' \bfC \bfz_{b} + 2\beta}}
	(\diag \bfA)^{-\frac{1}{2}} (\diag \bfB)^{\frac{1}{2}} \bfz_{b}.
\ea
Or, equivalently, 
\ba
	\left(\sqrt{\frac{a_{11}}{b_{11}}}\frac{\theta_1^{(y)}-\theta_{10}}{\sqrt{V_{11}^{(y)}}}, \ldots, 
	\sqrt{\frac{a_{mm}}{b_{mm}}}\frac{\theta_m^{(y)}-\theta_{m0}}{\sqrt{V_{mm}^{(y)}}}\right)'
	=\sqrt{\frac{m+2\alpha}{\bfz_{b}' \bfC \bfz_{b} + 2\beta}} ~\bfz_{b}.
\ea

The joint CDF of $H_1, \ldots, H_m$ is given by
 \ba
 	F(h_1, \ldots, h_m) &=& P(H_1 \leq h_1, \ldots, H_m \leq h_m) \\
	&=& P\left[ \Psi_{m+2\alpha}\left( \frac{\theta_1^{(y)}-\theta_{10}}{\sqrt{V_{11}^{(y)}}} \right) \leq h_1, 
		\ldots, \Psi_{m+2\alpha}\left( \frac{\theta_m^{(y)}-\theta_{m0}}{\sqrt{V_{mm}^{(y)}}} \right) \leq h_m \right] \\
	&=& P\left[ \frac{\theta_1^{(y)}-\theta_{10}}{\sqrt{V_{11}^{(y)}}} \leq \Psi_{m+2\alpha}^{-1}(h_1), 
		\ldots, \frac{\theta_m^{(y)}-\theta_{m0}}{\sqrt{V_{mm}^{(y)}}} \leq \Psi_{m+2\alpha}^{-1}(h_m) \right] \\
	&=& \Xi^m_{a, b}\left[ \sqrt{\frac{a_{11}}{b_{11}}} \Psi^{-1}_{m+2\alpha}(h_1), \cdots, 
		\sqrt{\frac{a_{mm}}{b_{mm}}} \Psi^{-1}_{m+2\alpha}(h_m) \right].
 \ea

If the misspecified covariance is used, follow the same argument in Lemma \ref{lemma1}.  
 \end{proof}

\subsection{Proof of Corollary \ref{coro1}}
\begin{proof}
Notice that $d\Phi^{-1}(x)/dx=1/\varphi(\Phi^{-1}(x))$, where $\varphi(x)=1/\sqrt{2\pi}e^{-x^2/2}$. Then
\ba
	f(h_i) &=& \varphi\left( \sqrt{r_i}\Phi^{-1}(h_i) \right) \frac{\sqrt{r_i}}{\varphi(\Phi^{-1}(h_i))} \\
	&=& \sqrt{r_i}\exp \left\{ \frac{1}{2}(1-r_i)[\Phi^{-1}(h_i)]^2 \right\}.
\ea
\end{proof}

\subsection{Proof of Corollary \ref{coro2}}
\begin{proof}
Let $u_i=\sqrt{r_i}\Phi^{-1}(h_i)$. Using the chain rule,
\ba
	f(h_1,\ldots,h_m) &=& \frac{\partial^m F}{\partial h_1 \cdots \partial h_m}=\frac{\partial^m \Phi_b^m}{\partial u_1 \cdots \partial u_m} 
	\prod_{i=1}^m \frac{\partial u_i}{\partial h_i} \\
	&=& \varphi_b^m(u_1, \ldots, u_m) \prod_{i=1}^m \frac{\sqrt{r_i}}{\varphi(\Phi^{-1}(h_i))} \\
	&=& \vert \bfP_b \vert^{-\frac{1}{2}}\exp\left\{ -\frac{1}{2}\bfu'\bfP_b^{-1}\bfu \right\} \prod_{i=1}^m \sqrt{r_i} \exp\left\{ \frac{1}{2} [\Phi^{-1}(h_i)]^2 \right\} \\
	&=& \left\{ \prod_{i=1}^m \sqrt{r_i} \right\} \vert \bfP_b \vert^{-\frac{1}{2}} \exp\left\{ -\frac{1}{2} \bfphi'\bfR^{\frac{1}{2}}\bfP_b^{-1} \bfR^{\frac{1}{2}} \bfphi 
	+\frac{1}{2} \bfphi'\bfphi \right\} \\
	&=& \left| \bfR^{\frac{1}{2}}\bfP_b^{-1} \bfR^{\frac{1}{2}} \right|^{\frac{1}{2}} 
	\exp \left\{ \frac{1}{2}\bfphi'(\bfI-\bfR^{\frac{1}{2}}\bfP_b^{-1} \bfR^{\frac{1}{2}})\bfphi \right\}. 
\ea
\end{proof}

\section*{Acknowledgement}
This research is partially supported by National Science Foundation under Grant No. OIA-1301789.

\renewcommand{\baselinestretch}{1.1}

%\small
%\bibliographystyle{agsm}
%\bibliography{ref}

\newpage

\begin{table}[h]
\caption{Soil Data: Disagreed decisions under Model 1 and 2, for the hypothesis $H_{0i}: \theta_i \geq 50$ versus $H_{0i}: \theta_i < 50$.
	$\dag$ denotes a rejection under the model. The nominal level is $0.05$.}
\begin{center}
\begin{tabular}{c|c|c|c}
\hline
\hline 
Time \& Location & Observed Value & Upper C.I. under Model 1 & Upper C.I. under Model 2 \\
\hline
January, Site 6 & 49.49 & 50.42$\dag$ & 51.03 \\
January, Site 70 & 49.55 & 50.52$\dag$ & 51.30 \\
January, Site 72 & 49.62 & 50.37$\dag$ & 50.89 \\
January, Site 87 & 49.70 & 50.67$\dag$ & 51.07 \\
January, Site 101 & 49.62 & 50.70$\dag$ & 51.21 \\
March, Site 8 & 49.94 & 50.78$\dag$ & 51.10 \\
March, Site 47 & 49.58 & 50.36$\dag$ & 51.22 \\
April, Site 38 & 49.53 & 50.19$\dag$ & 50.90 \\
June, Site 53 & 49.63 & 50.85$\dag$ & 51.83 \\
\hline
\end{tabular}
\end{center}
\label{tab:data}
\end{table}%

\newpage

\end{document}